%% file: non_displ_torus.tex
\documentclass[12pt,intlimits]{amsart}

\usepackage{hyperref}
\usepackage{graphicx}
\usepackage{color}
\usepackage{amsfonts,amssymb}
\usepackage{amsthm}
\usepackage{txfonts}

\usepackage{a4wide}

\newtheorem*{Cor*}{Corollary}

\newtheorem*{Thm*}{Theorem}

\newtheorem*{Prop*}{Proposition}
\theoremstyle{definition}

\newtheorem*{Rmk*}{Remark}

\newtheorem*{Ex*}{Example}

\newcommand{\Z}{\mathbb{Z}}
\newcommand{\R}{\mathbb{R}}

\newcommand{\RP}{\R\mathrm{P}}

\newcommand{\pf}{\longrightarrow}

\newcommand{\Mas}{\mu_{\mathrm{Maslov}}}

\renewcommand{\P}{\mathcal{P}}

\newcommand{\M}{\mathcal{M}}

\renewcommand{\H}{\mathrm{H}}

\newcommand{\HF}{\mathrm{HF}}

\newcommand{\Crit}{\mathrm{Crit}}

\newcommand{\beq}{\begin{equation}}
\newcommand{\beqn}{\begin{equation}\nonumber}
\newcommand{\eeq}{\end{equation}}

\newcommand{\bea}{\begin{equation}\begin{aligned}}
\newcommand{\bean}{\begin{equation}\begin{aligned}\nonumber}
\newcommand{\eea}{\end{aligned}\end{equation}}

\numberwithin{equation}{section}

\definecolor{Urs}{rgb}{1,0,0}
\definecolor{Peter}{rgb}{0,0,1}

\begin{document}
\title{A non-displaceable Lagrangian torus in $T^*S^2$}
\author{Peter Albers}
\author{Urs Frauenfelder}
\address{
    Peter Albers\\
    Courant Institute\\
    New York University}
\email{albers@cims.nyu.edu}
\address{
    Urs Frauenfelder\\
    Department of Mathematics\\
    LMU}
\email{Urs.Frauenfelder@mathematik.uni-muenchen.de}
\begin{abstract}
We show that the Lagrangian torus in the cotangent bundles of the 2-sphere obtained by applying the geodesic flow to the unit circle in
a fibre is not displaceable by computing its Lagrangian Floer homology. The computation is based on a symmetry argument.
\end{abstract}
\maketitle

Following a construction of Leonid Polterovich we consider the following Lagrangian torus $L$ in $T^*S^2$. Let $\varphi_t$ be the
geodesic flow on the cotangent bundle of the standard round $S^2$ where we identified $TS^2$ with $T^*S^2$ via the round metric. Fix
the unit circle $C\subset T^*_NS^2$ in the cotangent fiber over the north-pole $N\in S^2$. We consider the map
\bean
\phi:S^1\times S^1&\pf T^*S^2\\
(t,v)&\mapsto\varphi_t(v)
\eea
where we identify isometrically the circle $C\subset T^*_NS^2$ with $S^1$. The map $\phi$ defines a Lagrangian embedding, we set
\beqn
L:=\phi(T^2)\subset T^*S^2\;.
\eeq
Obviously, $\pi_2(T^*S^2,L)\cong\pi_2(T^*S^2)\oplus\pi_1(L)=\Z\oplus\Z\oplus\Z$, where the second $\Z$ is generated by the unit disk in $T_N^*S^2$
bounded by the loop $C$. The third disk corresponds to the loop $t\mapsto\varphi_t(v_0)$ in $\pi_1(L)$.

The disk bounded by $C$ lies in a Lagrangian submanifold namely $T_N^*S^2$ and thus has vanishing Maslov index and symplectic area. The
disk corresponding to the loop $t\mapsto\varphi_t(v_0)$ has Maslov index $2$.

We obtain a monotone Lagrangian torus $L$ with minimal Maslov number $N_L=2$. In particular, its Lagrangian Floer homology is
well-defined by standard means since $T^*S^2$ is a convex exact symplectic manifold and $L$ has minimal Maslov number $N_L\geq2$, see
\cite{Oh_Addendum_Floer_cohomology_of_Lagrangian_intersections_and_pseudo-holomorphic_disks_I}.

Apart from the fibers over the north-pole and the south-pole the torus $L$ intersects each fiber of $T^*S^2$ exactly twice. At north and south-pole
it intersects the fiber in a circle.\\[2ex]
The following question of Leonid Polterovich was posed to us by Felix Schlenk:

\begin{enumerate}
\item Is $L$ displaceable?
\item If not, what is $\HF(L,L)$?
\end{enumerate}

\begin{Thm*}The Floer homology of the torus $L$ is
\beqn
\HF_k(L,L)\cong\H_k(L;\Z/2\Z)\qquad k\in\Z/2\Z\,.
\eeq
In particular, $L$ is not displaceable in $T^*S^2$.
\end{Thm*}

\begin{Rmk*}
\begin{enumerate}
\item It was pointed out to us by L.~Polterovich that the non-displaceability of the Lagrangian torus gives another proof of the orderability
of the universal cover of the group of contactomorphism of $\RP^3$. Indeed, $T^*S^2$ can be regarded as the symplectization of $\RP^3$ and the
Lagrangian torus as pre-Lagrangian submanifold inside $\RP^3$. Thus, according to
\cite[Theorem 2.3.A]{Eliashberg_Polterovich_Partially_ordered_groups_and_geometry_of_contact_transformations} the orderability follows.
\item Recently, the following was proved by Paul Biran and Octav Cornea \cite{Biran_Cornea_Quantum_Structures}: The Lagrangian Floer homology
of a monotone Lagrangian torus either vanishes identically or is isomorphic to the singular homology of the torus.
\end{enumerate}
\end{Rmk*}

\begin{proof}
We consider the symplectic involution $I$ on $T^*S^2$ which is induced by the reflection across a great circle through the north-pole.
Since this reflection is an isometry of $S^2$ the involution $I$ leaves $L$ invariant: $I(L)=L$.

On the other hand the fixed point set of $I$ can be identified with $T^*S^1$. Furthermore, $L\cap T^*S^1$ consists out of two copies of $S^1$
which are symmetric under reflection at the zero-section of $T^*S^1$.

We choose on $S^1=\R/\Z$ the height function $h(t):=\cos(2\pi t)$. On the circle $C\subset T^*_NS^2$ we fix the directions tangential to
the above chosen great circle to be the critical points of $h$ and set $f(t,v):=\varepsilon\big(h(t)+h(v)\big)$ where $\varepsilon>0$ is chosen below.

Note that $f$ is a Morse function on $L$ which is invariant under the symplectic involution $I$. Furthermore, the critical points of $f$ are
fixed under the action of $I$, see figure \ref{fig:sphere}.

\begin{figure}[ht]
\input{sphere.tex}
\caption{}\label{fig:sphere}
\end{figure}

From the proof of Weinstein's Lagrangian neighborhood theorem in \cite[Theorem 3.33]{McDuff_Salamon_introduction_symplectic_topology} it is
apparent that the neighborhood can be chosen $I$-invariant. (Simply choose the compatible almost complex structure appearing in the proof to
be $I$-invariant.)
Then we extend $f$ independently of the fiber-coordinate of the neighborhood to a compactly supported $I$-invariant function on $T^*S^2$, which we
call $f$ again.

If we choose $\varepsilon$ sufficiently small, the set
\beqn
\P_L(f):=\{x\in C^\infty([0,1],T^*S^2)\mid \dot{x}=X_f(x),\,x(0),x(1)\in L\}
\eeq
equals $\Crit(f|_L)$ and thus is point-wise fixed by $I$.

We choose a time-dependent family $J_t$ of (compatible) almost complex structures on $T^*S^2$ and consider
perturbed holomorphic strips, i.e.~smooth maps $u:\R\times[0,1]\pf M$ solving
\beqn
\tag{$\star$}\partial_s u +J_t(u)\big(\partial_tu-X_f(u)\big)=0,\quad u(s,0),\,u(s,1)\in L\quad\text{and}\quad u(\pm\infty)=x_\pm\in\Crit(f|_L)\;.
\eeq
For regular almost complex structures the set $\M(x_-,x_+;J)$ of these maps is a smooth finite dimensional manifold.
Note that the Hamiltonian vector field $X_f$ is $I$-invariant. If we choose the family of almost complex structures to be $I$-invariant (see below)
we obtain an induced action on the moduli-spaces $\M(x_-,x_+;J)$.
The fixed points of this action are those perturbed strips $u$ such that $u(\R\times[0,1])\subset\text{Fix}(I)\cong T^*S^1$ and $u$ has boundary
on $L\cap\text{Fix}(I)=S^1\sqcup S^1$.

The moduli spaces $\M(x_-,x_+;J)$ consist out of many connected components in general, which we coarsely separate by the Maslov index of
the relative homotopy class of the perturbed strip
\beqn
\M(x_-,x_+;J,n):=\{u\in\M(x_-,x_+;J)\mid \Mas([u])=n\}\;.
\eeq
In particular, we have $\M(x_-,x_+;J)=\bigcup_n\M(x_-,x_+;J,n)$ and $\dim \M(x_-,x_+;J,n)=i(x_-)-i(x_+)+n$, where
$i(x_-)$ is the Morse index of the critical point $x_-$ and
where we implicitly assume that $\varepsilon$ is chosen small enough.

By definition the boundary operator $\partial_F$ in Floer theory counts (in $\Z/2\Z$) the number of elements of $\M(x_-,x_+;J,n)/\R$ in case
$\dim\M(x_-,x_+;J,n)=1$.

In \cite[Proposition 5.13]{Khovanov_Seidel_quivers_Floer_cohomology_and_braid_group_actions} it is proved that the set of $I$-invariant families of
almost complex structures which are regular \textit{outside} the fixed point set of $I$ is a generic subset of the set of $I$-invariant families of
almost complex structures. We choose such a family of almost complex structures $\widehat{J}$.
If $u$ is fixed by the $I$-action it satisfies $u(\R\times[0,1])\subset\text{Fix}(I)\cong T^*S^1$ with boundary on $L\cap\text{Fix}(I)=S^1\sqcup S^1$
and $u(\pm\infty)=x_\pm$. Since $\pi_2(T^*S^1,S^1\sqcup S^1)=0$ the map $u$ is contractible relative $L$ in $T^*S^1$ to a point and hence also in the
ambient space. In particular, $\Mas(u)=0$. We conclude that for $n\not=0$ no element $u\in\M(x_-,x_+;\widehat{J},n)$ can lie in the fixed point
set of $I$.

Since the family of almost complex structures is regular outside the fixed point set we obtain smooth moduli spaces $\M(x_-,x_+;\widehat{J},n)$ with
a free $I$-action. In particular, if $\dim\M(x_-,x_+;\widehat{J},n)=1$ the space $\M(x_-,x_+;\widehat{J},n)/\R$ is compact and contains an even
number of elements.\footnote{L.~Polterovich pointed out that fixed points of $I$ on the quotient $\M(x_-,x_+;\widehat{J},n)/\R$ are equivalence
classes $\{u\}$ of maps satisfying $I(u(s+\sigma,t))=u(s,t)$. But this implies that $u$ is actually periodic in $s$ which contradicts the
convergence at $\pm\infty$.}

The $I$-action on the torus $L$ is depicted in figure \ref{fig:torus}. The dots are the four critical points on the torus which is obtained by identifying
opposite sides of the square. The critical points from figure \ref{fig:sphere} correspond to the dots on the center horizontal lines.

\begin{figure}[ht]
\input{torus.tex}
\caption{}\label{fig:torus}
\end{figure}

We observe that miraculously the standard flat metric on $T^2$ is Morse-Smale for $f$ \textit{and} $I$-invariant. In general, there is no reason
to expect the existence of an invariant metric which is Morse-Smale. The metric determines a canonical almost complex structure $J_{LC}$ on $T^*L$
which is called Levi-Civita almost complex structure. Since we chose the Weinstein neighborhood to be $I$-invariant we can extend $J_{LC}$ to
an $I$-invariant almost complex structure on $T^*S^2$ which we again denote by $J_{LC}$.

By an argument of Po\'zniak \cite{Pozniak_PhD} the number of elements in $\M(x_-,x_+;J_{LC},0)/\R$ is given by the number of Morse trajectories connecting
$x_-$ and $x_+$. Indeed, because in the Weinstein neighborhood of $L$ the Hamiltonian function $f$ is constant along the fibers it follows for small enough $\varepsilon$
from Po\'zniaks theorem \cite{Pozniak_PhD} that all elements $u\in\M(x_-,x_+;J_{LC},0)$ equal Morse flow lines and thus $J_{LC}$ is regular
for all $u\in\M(x_-,x_+;J_{LC},0)$. By \cite[Proposition 5.13]{Khovanov_Seidel_quivers_Floer_cohomology_and_braid_group_actions} we can
find an $I$-invariant almost complex structure $\widehat{J}$ arbitrarily close to $J_{LC}$ which is regular for all
$u\in\M(x_-,x_+;\widehat{J},n\not=0)$. Furthermore, using the
implicit function theorem we can canonically identify $\M(x_-,x_+;\widehat{J},0)$ with $\M(x_-,x_+;J_{LC},0)$, hence $\widehat{J}$ is regular
throughout.

Recall the definition of the Floer boundary operator $\partial_F$ on an element $x_-\in\P_L(f)=\Crit(f|_L)$
\beqn
\partial_F x_-=\sum_{x_+}\#_2\M(x_-,x_+;\widehat{J})/\R\cdot x_+\;,
\eeq
where $\#_2\M(x_-,x_+;\widehat{J})/\R$ is the number in $\Z/2\Z$ of elements in the zero-dimensional components. Equivalently,
\bean
\partial_F x_-&=\sum_{x_+,\;n}\#_2\M(x_-,x_+;\widehat{J},n)/\R\cdot x_+\\
            &=\sum_{x_+,\;n\not=0}\#_2\M(x_-,x_+;\widehat{J},n)/\R\cdot x_++\sum_{x_+}\#_2\M(x_-,x_+;\widehat{J},0)/\R\cdot x_+\;.
\eea
The second summand equals the Morse differential $\partial_M$ and the first summand vanishes in $\Z/2\Z$ by the above discussion.
In particular, we obtain
\beqn
\partial_F=\partial_M\;.
\eeq
This proves the theorem.
\end{proof}

\begin{Rmk*}
In higher dimensions a Lagrangian $L:=S^1\times S^{n-1}$ can be constructed by applying the geodesic flow to the unit sphere in $T^*_NS^n$.
One obvious generalization of the above approach is to consider a reflection on a hyper-sphere. This indeed defines a symplectic involution
$I_n$ on $T^*S^n$ but with $\mathrm{Fix}(I_n)=T^*S^{n-1}$. Thus, there exist non-trivial perturbed holomorphic strips mapping into the fixed point set
and the above approach fails.

The following was remarked by Florin Belgun and Slava Matveev. Let us again fix the sum of the height functions on $S^1\times S^{n-1}$. This
defines a Morse function on $L$ where we require the critical points to lie as antipodal points in the fibres over the north and south pole. Thus,
the two critical points in $T^*_NS^n$ define a unique geodesic $S^1$ through north and south pole. Denote by $I$ the reflection on this $S^1$.
This induces a symplectic involution on $T^*S^n$ leaving $L$ invariant and fixing the critical points. Furthermore, the fixed point set is
$T^*S^1$. Therefore, the above proof can be repeated verbatim.
\end{Rmk*}

\subsubsection*{Acknowledgments}
This result was proved during the stay of the second author at the Courant Institute. Both authors thank the Courant Institute for its stimulating
working atmosphere.

Both author are partially supported by the German Research Foundation (DFG) through Priority Programm 1154
"Global Differential Geometry", grants AL 904/1-1 and CI 45/1-2 and by NSF Grants DMS-0102298 and DMS-0603957.

\noindent\hrulefill
%
%
\bibliographystyle{amsalpha}
\bibliography{../../Bibtex/bibtex_paper_list}
\end{document}

%% file: sphere.tex
\begin{picture}(0,0)%
\includegraphics{sphere.pstex}%
\end{picture}%
\setlength{\unitlength}{3947sp}%
\begingroup\makeatletter\ifx\SetFigFont\undefined%
\gdef\SetFigFont#1#2#3#4#5{%
  \reset@font\fontsize{#1}{#2pt}%
  \fontfamily{#3}\fontseries{#4}\fontshape{#5}%
  \selectfont}%
\fi\endgroup%
\begin{picture}(5124,4431)(3451,-5023)
\put(5976,-1351){\makebox(0,0)[lb]{\smash{{\SetFigFont{12}{14.4}{\rmdefault}{\mddefault}{\updefault}{\color[rgb]{0,0,0}$C$}%
}}}}
\put(5701,-3511){\makebox(0,0)[lb]{\smash{{\SetFigFont{12}{14.4}{\rmdefault}{\mddefault}{\updefault}{\color[rgb]{0,0,0}$I$}%
}}}}
\put(7576,-1336){\makebox(0,0)[lb]{\smash{{\SetFigFont{12}{14.4}{\rmdefault}{\mddefault}{\updefault}{\color[rgb]{0,0,0}$T^*_NS^2$}%
}}}}
\put(5626,-736){\makebox(0,0)[lb]{\smash{{\SetFigFont{12}{14.4}{\rmdefault}{\mddefault}{\updefault}{\color[rgb]{0,0,0}critical point}%
}}}}
\put(3451,-1936){\makebox(0,0)[lb]{\smash{{\SetFigFont{12}{14.4}{\rmdefault}{\mddefault}{\updefault}{\color[rgb]{0,0,0}critical point}%
}}}}
\end{picture}%

%% file: torus.tex
\begin{picture}(0,0)%
\includegraphics{torus.pstex}%
\end{picture}%
\setlength{\unitlength}{3947sp}%
\begingroup\makeatletter\ifx\SetFigFont\undefined%
\gdef\SetFigFont#1#2#3#4#5{%
  \reset@font\fontsize{#1}{#2pt}%
  \fontfamily{#3}\fontseries{#4}\fontshape{#5}%
  \selectfont}%
\fi\endgroup%
\begin{picture}(1974,3106)(4414,-4205)
\put(5318,-4156){\makebox(0,0)[lb]{\smash{{\SetFigFont{12}{14.4}{\rmdefault}{\mddefault}{\updefault}{\color[rgb]{0,0,0}$I$}%
}}}}
\end{picture}%

%% file: non_displ_torus.bbl
\providecommand{\bysame}{\leavevmode\hbox to3em{\hrulefill}\thinspace}
\providecommand{\MR}{\relax\ifhmode\unskip\space\fi MR }
\providecommand{\MRhref}[2]{%
  \href{http://www.ams.org/mathscinet-getitem?mr=#1}{#2}
}
\providecommand{\href}[2]{#2}
\begin{thebibliography}{Po{\'z}99}

\bibitem[BC]{Biran_Cornea_Quantum_Structures}
P.~Biran and O.~Cornea, in preparation.

\bibitem[EP00]{Eliashberg_Polterovich_Partially_ordered_groups_and_geometry_of%
_contact_transformations}
Y.~Eliashberg and L.~Polterovich, \emph{Partially ordered groups and geometry
  of contact transformations}, Geom. Funct. Anal. \textbf{10} (2000), no.~6,
  1448--1476.

\bibitem[KS02]{Khovanov_Seidel_quivers_Floer_cohomology_and_braid_group_action%
s}
M.~Khovanov and P.~Seidel, \emph{Quivers, {F}loer cohomology, and braid group
  actions}, J. Amer. Math. Soc. \textbf{15} (2002), no.~1, 203--271
  (electronic).

\bibitem[MS98]{McDuff_Salamon_introduction_symplectic_topology}
D.~McDuff and D.~Salamon, \emph{Introduction to symplectic topology}, second
  ed., Oxford Mathematical Monographs, The Clarendon Press Oxford University
  Press, New York, 1998.

\bibitem[Oh95]{Oh_Addendum_Floer_cohomology_of_Lagrangian_intersections_and_ps%
eudo-holomorphic_disks_I}
Y.-G. Oh, \emph{Addendum to: ``{F}loer cohomology of {L}agrangian intersections
  and pseudo-holomorphic disks {I}.''}, Comm. Pure Appl. Math. \textbf{48}
  (1995), no.~11, 1299--1302.

\bibitem[Po{\'z}99]{Pozniak_PhD}
M.~Po{\'z}niak, \emph{Floer homology, {N}ovikov rings and clean intersections},
  Northern California Symplectic Geometry Seminar, Amer. Math. Soc. Transl.
  Ser. 2, vol. 196, Amer. Math. Soc., Providence, RI, 1999, pp.~119--181.

\end{thebibliography}
